\documentstyle[leqno,12pt]{article}

\title{Nonlinear Schr\"odinger equations with\\
radially symmetric data of critical regularity}

\author{Kunio Hidano}
\date{}

\begin{document}
\maketitle

\begin{abstract}
This paper is concerned with 
the global existence of small solutions 
to pure-power nonlinear Schr\"odinger equations subject to 
radially symmetric data with critical regularity. 
Under radial symmetry 
we focus our attention 
on the case where the power of nonlinearity is somewhat smaller 
than the pseudoconformal power 
and the initial data belong to the 
scale-invariant homogeneous Sobolev space. 
In spite of the negative-order differentiability of initial data 
the nonlinear Schr\"odinger equation 
has global in time solutions 
provided that 
the initial data have the small norm. 
The key ingredient in the proof of this result is an effective use of 
global weighted smoothing estimates 
specific to radially symmetric solutions.


{\it {\rm 2000} Mathematics Subject Classification Numbers.} 35Q55, 35B65.
\end{abstract}
\section{Introduction and result}
\baselineskip=0.55cm
In this paper we study the Cauchy problem 
for the nonlinear Schr\"odinger equation
\begin{equation}
i\partial_t u+\Delta u
=
\lambda |u|^{p-1}u,
\quad
t\in{\bf R},\,\,x\in{\bf R}^n
\end{equation}
$(p>1,\,\,\lambda\in{\bf C})$ 
subject to the initial data $u(0,x)=\varphi(x)$. 
Our main concern is to investigate the problem of 
global existence of small $H^s$-solutions. 
Assuming
\begin{eqnarray}
& &
\hspace{2cm}
p-1=\frac{4}{n-2s},\,\,\,
0\leq s<\frac{n}{2}\\
& &
\biggl(
\mbox{or equivalently}\,\,
s=\frac{n}{2}-\frac{2}{p-1},\,\,\,\,
1+\frac{4}{n}\leq p<\infty
\biggr)\nonumber
\end{eqnarray}
and in addition 
$[s]<p-1$ 
($[s]$ denotes the greatest integer 
not greater than $s\geq 0$) 
if $p-1$ is not an even integer, 
Cazenave and Weissler proved that, 
for any $\varphi\in H^s({\bf R}^n)$ with 
$\||D_x|^s\varphi\|_{L^2({\bf R}^n)}$ small, 
the associated integral equation 
has a unique global solution 
(Theorems 1.1 and 1.2 in \cite{CW}). 
From the point of view of 
the $H^s$-theory, 
their result concerns the critical case. 
Though they left it open to show 
the global existence 
of small $H^s$-solutions 
in the case where subcritical nonlinear terms 
occur in the equation 
together with the critical nonlinear term 
(see page 81 of \cite{CW}), 
this issue was resolved by Kato 
to some extent (Theorem 6.1 in \cite{K}). 
Note that, 
for the proof of small data global existence, 
Kato assumed much less than Cazenave and Weissler did, 
but still he assumed the restrictive condition (F.3) 
on page 299 of \cite{K} 
which keeps any $L^2$-subcritical term $\lambda|u|^{p-1}u$ 
with $p<1+4/n$ from occurring in the equation. 
We also refer to the fine works of Ginibre, Ozawa and Velo \cite{GOV}, 
Nakamura and Ozawa \cite{NO1}, and Pecher \cite{P2}. 
In these papers they also devised how to 
relax the assumption on the nonlinear term 
Cazenave and Weissler made in \cite{CW}. 
Allowing the $L^2$-critical and/or the $L^2$-supercritical terms 
in the equation, 
they also established theorems on the global existence of small solutions 
in the $H^s$-framework for $0\leq s<n/2$. 
For the $H^s$-theory with $s\geq n/2$ we refer the reader to 
Nakamura and Ozawa \cite{NO98}, \cite{NO2000}.

In this paper we focus our attention 
on radially symmetric solutions 
to the pure-power nonlinear Schr\"odinger equation, 
and we intend to explore the subject concerning 
the global existence of small solutions 
when $p$ is somewhat smaller than the $L^2$-critical power $1+4/n$ 
and the radially symmetric initila data $\varphi$ is in the scale-invariant 
homogeneous Sobolev space. 
Let $U(t)=e^{it\Delta}$, 
and let ${\dot H}^s_2({\bf R}^n)$ denote the homogeneous Sobolev space 
$\{\,v\in{\mathcal S}'({\bf R}^n)\,|\,|D_x|^sv\in L^2({\bf R}^n)\,\}$ 
for $-n/2<s<n/2$, 
where 
$|D_x|^sv:={\mathcal F}^{-1}|\xi|^s{\mathcal F}v$. 
We shall prove

\vspace{0.3cm}

\noindent{\bf Theorem 1.1.} {\it 
Suppose that $n\geq 3$, $4/(n+1)<p-1<4/n$, $\lambda\in{\bf C}$. 
Set $s_0=-n/2+2/(p-1)$. 
There exists a positive constant $\delta$ depending on 
$n$, $p$, $\lambda$ such that 
if radially symmetric data 
$\varphi\in{\dot H}^{-s_0}_2({\bf R}^n)$ 
is small so that 
$\||D_x|^{-s_0}\varphi\|_{L^2({\bf R}^n)}\leq \delta$, then 
the integral equation
\begin{equation}
u(t)
=
U(t)\varphi
-i\lambda
\int_0^t
U(t-\tau)
|u(\tau)|^{p-1}u(\tau)d\tau
\end{equation}
has a unique radially symmetric solution 
$u\in C({\bf R};{\dot H}^{-s_0}_2({\bf R}^n))$ 
satisfying 
\begin{equation}
\sup_{t\in{\bf R}}
\||D_x|^{-s_0}u(t,\cdot)\|_{L^2({\bf R}^n)}
+
\||x|^{-\alpha_0}u\|_{L^{q_0}({\bf R}\times{\bf R}^n)}
\leq 
C\delta
\end{equation}
for a suitable constant $C>0$. 
Here $q_0>2$ and $\alpha_0\in{\bf R}$ are the exponents defined 
in $(3.5)$, $(3.6)$ below.
}

\vspace{0.3cm}

Note that $0<s_0<1/2$ for $4/(n+1)<p-1<4/n$. 
Hence the key to showing this theorem is 
getting over the difficulty caused by 
the negative-order differentiability of initial data. 
The Strichartz estimates manifest 
the gain of {\it integrability exponents} 
for the free solutions with initial data in $L^2({\bf R}^n)$ \cite{Str}. 
They have played an essential role 
in the development of the local and global (in time) $H^s$-theory 
for $s\geq 0$, together with 
space-time estimates of solutions to inhomogeneous equations 
as well as elaborate estimates of products of functions 
in fractional-order Sobolev or Besov spaces. 
For the purpose of carrying out the contraction-mapping argument 
in the present setting 
we therefore start by showing the Strichartz-type estimates 
with {\it derivative gain} for the free Schr\"odinger equation 
subject to initial data in $L^2({\bf R}^n)$ (see Section 2). 
In order to obtain such estimates 
we adapt a technique of showing extended Strichartz estimates 
for the free wave equation 
with radially symmetric data 
(see, e.g., Sogge \cite{So} on page 126), 
and in this way we shall prove weighted space-time 
$L^q({\bf R}^{1+n})$ estimates 
with derivative gain for solutions 
to the free Schr\"odinger equation 
with radially symmetric data in $L^2({\bf R}^n)$. 
In application to nonlinear problems 
it is also necessary to establish space-time estimates 
for the inhomogeneous Schr\"odinger equation. 
By virtue of the Christ-Kiselev lemma \cite{CK} 
along with the $TT^*$ argument, 
it is possible to derive 
some weighted estimates 
of radially symmetric solutions 
to the inhomogeneous equation 
directly from the weighted estimates 
for the free Schr\"odinger equation. 
Making use of these weighted estimates, 
we can carry out the contraction-mapping 
argument to show the main result. 

In connection with Theorem 1.1 the referee has advised the author 
to make reference to the results 
of Y.\,Tsutsumi \cite{Tsu}, 
Ginibre, Ozawa and Velo \cite{GOV}, 
and Nakanishi and Ozawa \cite{NaOz} 
who also studied the global existence for the nonlinear 
Schr\"odinger equation with the $L^2$-subcritical term 
$\lambda |u|^{p-1}u$ $(0<p-1<4/n,\,\lambda\in{\bf C})$. 
For $\lambda\in{\bf R}$ it was proved that 
the nonlinear Schr\"odinger equation admits a unique global solution 
for any data $\varphi\in L^2({\bf R}^n)$ \cite{Tsu}. 
When $\lambda\in{\bf C}$ and 
$$
\frac2n<p-1<\frac4n\,\,\,(n\leq 3),\,\,\,\,\,
\frac{\sqrt{n^2+4n+36}-n-2}{4}<p-1<\frac4n\,\,\,(n\geq 4),
$$
it follows from Theorem 2.1 of \cite{NaOz}, 
which refines an earlier result of \cite{GOV}, 
that the nonlinear Schr\"odinger equation 
admits a unique global solution 
for any data $\varphi\in L^2({\bf R}^n)\cap{\cal F}{\dot H}^{s_0}_2$ 
having its small ${\cal F}{\dot H}^{s_0}_2$-norm. 
Here $s_0$ is the same as in Theorem 1.1. 
In view of the dual Hardy inequality 
$$
\||D_x|^{-s_0}v\|_{L^2({\bf R}^n)}
\leq
C
\||x|^{s_0}v\|_{L^2({\bf R}^n)}
$$
we find that under the spherical symmetry assumption 
Theorem 1.1 improves on the result of Nakanishi and Ozawa 
because it uses the weaker norm to measure the size of initial data. 
The author thanks the referee 
for drawing attention to the relation 
between Theorem 1.1 and Theorem 2.1 of \cite{NaOz}. 

Finally we mention a similar progress on 
the $H^s$-theory for the nonlinear wave equation
\begin{equation}
\partial_t^2 u-\Delta u
=
\lambda|u|^{p-1}u.
\end{equation}
Global existence of small $H^s$-solutions 
has been studied in Lindblad and Sogge \cite{LS}, 
Nakamura and Ozawa \cite{NO2}. 
(See also \cite{LS}, \cite{Kap}, \cite{P1}, \cite{NO2} on 
the local in time $H^s$-theory.) 
Some improvements have been made in the radial case 
thanks to extended Strichartz-type estimates 
specific to radially symmetric solutions. 
See Theorem 6.6.1 of Sogge \cite{So} 
and Hidano \cite{H1}, \cite{H2} for this matter. 

This paper is organized as follows. 
Section 2 is devoted to the proof of weighted space-time 
$L^q$ estimates for the free Schr\"odinger equation and 
the inhomogeneous equation. 
In Section 3 we prove Theorem 1.1 by using these key estimates. 
\section{Weighted space-time $L^q$ estimates}
\setcounter{equation}{0}
In this section we first prove weighted $L^q$ estimates 
for the free Schr\"odinger equation with radially symmetric data 
by adapting the method of showing improved radial Strichartz estimates 
for the free wave equation (see, e.g., Sogge \cite{So} on page 126). 
Subsequently, using the celebrated lemma of Christ and Kiselev \cite{CK} 
(see also Smith and Sogge \cite{SS}, Tao \cite{T}) 
along with the $TT^*$ argument, 
we derive some weighted $L^{{\tilde q}'}$-$L^q$ estimates 
for the inhomogeneous equation 
from the weighted estimate for the free Schr\"odinger equation 
with radially symmetric data. 

Let us start with showing the following theorem.

\vspace{0.3cm}

\noindent {\bf Theorem 2.1.} {\it Suppose $n\geq 2$. 
There exists a constant $C>0$ depending on 
$n$, $q$, $\alpha$, and the estimate 
\begin{equation}
\||x|^{-\alpha}|D_x|^sU(t)\varphi\|
_{L^q({\bf R}\times{\bf R}^n)}
\leq
C\|\varphi\|_{L^2({\bf R}^n)}
\end{equation}
holds for radially symmetric data $\varphi\in L^2({\bf R}^n)$ 
provided that 
\begin{equation}
-\alpha-s+\frac{n+2}{q}=\frac{n}{2},
\,\,\,
2\leq q<\infty,\,\,\,
\frac{n}{q}-\frac{n-1}{2}<\alpha<\frac{n}{q}.
\end{equation}
}

\noindent{\it Proof}. 
We follow the proof of Strichartz-type estimates 
for the free wave equation with radially symmetric data. 
(See page 126 of \cite{So}. See also \cite{H1}, \cite{H2}.) 
Let $d\sigma$ denote the Lebesgue measure on the unit sphere 
$S^{n-1}\subset{\bf R}^n$, and 
let $\widehat{d\sigma}$ denote the Fourier transform of $d\sigma$, that is, 
\begin{equation}
\widehat{d\sigma}(\xi)
=
\int_{S^{n-1}_\omega}e^{-i\omega\cdot\xi}d\sigma
\quad
(d\sigma=d\sigma(\omega)).
\end{equation}
Since $\varphi$ is radially symmetric 
and so is its Fourier transform $\hat\varphi$, 
we may write $\hat\varphi(\xi)$ as $\psi(|\xi|)$. 
Making the change of variables $\eta=\rho^2$, 
we have 
\begin{eqnarray}
& &
|D_x|^s(U(-t)\varphi)(x)\\
& &
=(2\pi)^{-n}
\int_0^\infty
e^{it\rho^2}\psi(\rho)\rho^{s+n-1}
\widehat{d\sigma}(\rho x)d\rho\nonumber\\
& &
=(2\pi)^{-n}
\int_0^\infty
e^{it\eta}\psi(\sqrt{\eta})\eta^{(s+n-1)/2}
\widehat{d\sigma}(\sqrt{\eta} x)
\frac{1}{2\sqrt{\eta}}d\eta\nonumber\\
& &
=(2\pi)^{-n+1}
\biggl(
(2\pi)^{-1}
\int_{\bf R}
e^{it\eta}H(\eta)\psi(\sqrt{|\eta|})|\eta|^{(s+n-1)/2}
\widehat{d\sigma}(\sqrt{|\eta|} x)
\frac{1}{2\sqrt{|\eta|}}d\eta
\biggr).\nonumber
\end{eqnarray}
Here $H(\eta)$ denotes the Heaviside function. 
We therefore get by the Sobolev embedding theorem and 
the Plancherel theorem 
\begin{eqnarray}
& &
\||D_x|^s(U(\cdot)\varphi)(x)\|_{L^q({\bf R}_t)}\\
& &
\leq
C
\biggl\|
|D_t|^{1/2-1/q}\nonumber\\
& &
\hspace{2cm}
\biggl(
(2\pi)^{-1}
\int_{\bf R}
e^{it\eta}H(\eta)\psi(\sqrt{|\eta|})|\eta|^{(s+n-2)/2}
\widehat{d\sigma}(\sqrt{|\eta|} x)
d\eta
\biggr)
\biggr\|_{L^2({\bf R}_t)}\nonumber\\
& &
=
C\|\eta^{1/2-1/q+(s+n-2)/2}\psi(\sqrt\eta)
\widehat{d\sigma}(\sqrt{\eta} x)\|_{L^2({\bf R}_\eta^+)}\nonumber
\end{eqnarray}
$({\bf R}^+=(0,\infty))$ for any fixed $x\in{\bf R}^n$. 
Using the Minkowski inequality and 
making the change of variables 
$y:=\sqrt\eta x$, 
we obtain from (2.5)
\begin{eqnarray}
& &
\||x|^{-\alpha}|D_x|^s U(t)\varphi\|_{L^q({\bf R}_t\times{\bf R}_x^n)}\\
& &
\leq
C
\bigl\|
\eta^{1/2-1/q+(s+n-2)/2}\psi(\sqrt\eta)
\||x|^{-\alpha}\widehat{d\sigma}(\sqrt\eta x)\|
_{L^q({\bf R}^n_x)}
\bigr\|_{L^2({\bf R}^+_\eta)}\nonumber\\
& &
=C
\bigl\|
\eta^{1/2-1/q+(s+n-2)/2}\psi(\sqrt\eta)\eta^{\alpha/2-n/(2q)}
\||y|^{-\alpha}\widehat{d\sigma}(y)\|_{L^q({\bf R}^n_y)}
\bigr\|_{L^2({\bf R}^+_\eta)}\nonumber\\
& &
\leq
C\|\eta^{1/2-1/q+(s+n-2)/2+\alpha/2-n/(2q)}
\psi(\sqrt\eta)\|_{L^2({\bf R}^+_\eta)}.\nonumber
\end{eqnarray}
Here we have handled the $L^q({\bf R}^n_y)$-norm of 
$|y|^{-\alpha}\widehat{d\sigma}(y)$ as 
\begin{equation}
\||y|^{-\alpha}\widehat{d\sigma}(y)\|
_{L^q(\{y\in{\bf R}^n:|y|<1\})}
+
\||y|^{-\alpha}\widehat{d\sigma}(y)\|
_{L^q(\{y\in{\bf R}^n:|y|>1\})}
\leq
C,
\end{equation}
using the assumption 
$(n/q)-(n-1)/2<\alpha<n/q$ and the fact that 
$\widehat{d\sigma}$ is a smooth function 
on ${\bf R}^n$ satisfying 
$\widehat{d\sigma}(y)=O(|y|^{-(n-1)/2})$ 
as $|y|\to+\infty$ 
(See, e.g., Stein \cite{Ste} on page 348.). 
To finish the proof,
we make the change of variables 
$\lambda=\sqrt\eta$, 
and we continue the estimate of (2.6) as 
\begin{eqnarray}
& &
\dots
=C\|\lambda^{1-(2/q)+s+n-2+\alpha-(n/q)}\psi(\lambda)\lambda^{1/2}\|
_{L^2({\bf R}^+_\lambda)}
\\
& &
\hspace{0.6cm}
=C\|\lambda^{(n-1)/2}\psi(\lambda)\|_{L^2({\bf R}^+_\lambda)}
=
C\|\varphi\|_{L^2({\bf R}^n)}.\nonumber
\end{eqnarray}
We have finished the proof of Theorem 2.1.
$\hfill\Box$

\vspace{0.2cm}

\noindent{\it Remark.} The referee has kindly pointed out that 
Theorem 2.1 is true also for $q=\infty$. 
Indeed, using the weighted Sobolev inequality for 
radially symmetric functions (see, e.g., Appendix of \cite{H0}), 
we get for all $(t,x)\in{\bf R}\times{\bf R}^n$
\begin{eqnarray*}
& &
|x|^{(n/2)-s}|(U\varphi)(t,x)|\\
& &
\leq
C\||D_x|^s(U\varphi)(t,\cdot)\|_{L^2({\bf R}^n)}
=
C\||D_x|^s\varphi\|_{L^2({\bf R}^n)},
\quad
\frac12<s<\frac{n}{2}.
\end{eqnarray*}
This immediately leads to the estimate (2.1) for $q=\infty$. 
The referee has also pointed out that 
it is hence possible to derive the estimate (2.1) for $2<q<\infty$ 
via the complex interpolation between 
the above case $q=\infty$ and the well-known estimate 
for $q=2$ obtained in $\cite{KY}$, $\cite{BA-K}$, $\cite{Su}$, and $\cite{V}$. 

\vspace{0.2cm}

Using the Christ-Kiselev lemma \cite{CK} 
(see also Smith and Sogge \cite{SS}, Tao \cite{T}), 
we can derive the following estimates 
for the inhomogeneous equation 
from the estimate (2.1) for the free solution $U(t)\varphi$. 

\vspace{0.3cm}

\noindent{\bf Theorem 2.2.} {\it Suppose 
$n\geq 2$, $2\leq q<\infty$, $2\leq{\tilde q}<\infty$, 
$n/q-(n-1)/2<\alpha<n/q$, 
$n/{\tilde q}-(n-1)/2<{\tilde\alpha}<n/{\tilde q}$. 
Set
$$
s:=-\alpha+\frac{n+2}{q}-\frac{n}{2},\,\,\,
{\tilde s}:=-\tilde\alpha+\frac{n+2}{\tilde q}-\frac{n}{2}.
$$
The estimate
\begin{equation}
\biggl\|
|x|^{-\alpha}|D_x|^s
\int_0^tU(t-\tau)F(\tau)d\tau
\biggr\|
_{L^q({\bf R}\times{\bf R}^n)}
\leq
C\|
|x|^{\tilde\alpha}|D_x|^{-\tilde s}F\|
_{L^{{\tilde q}'}({\bf R}\times{\bf R}^n)}
\end{equation}
holds for radially symmetric 
$($in $x$$)$ $F$.
}

\vspace{0.3cm}

\noindent{\bf Theorem 2.3.} {\it 
Suppose $n\geq 2$, $2\leq{\hat q}<\infty$, 
$n/{\hat q}-(n-1)/2<{\hat \alpha}<n/{\hat q}$. 
Set
${\hat s}:=-\hat\alpha+(n+2)/{\hat q}-n/2$. 
The estimate 
\begin{equation}
\biggl\|
|D_x|^{\hat s}
\int_0^tU(t-\tau)F(\tau)d\tau
\biggr\|
_{L^\infty({\bf R};L^2({\bf R}^n))}
\leq
C\|
|x|^{\hat\alpha}F\|
_{L^{{\hat q}'}({\bf R}\times{\bf R}^n)}
\end{equation}
holds for radially symmetric 
$($in $x$$)$ $F$.
}

\vspace{0.3cm}

\noindent{\bf Corollary 2.4.} {\it 
Suppose $n\geq 2$, 
$2\leq q<\infty$, 
$2\leq{\tilde q}<\infty$, 
$n/q-(n-1)/2<\alpha<n/q$, 
$n/{\tilde q}-(n-1)/2<{\tilde\alpha}<n/{\tilde q}$, and 
$$
-\alpha+\frac{n+2}{q}-\frac{n}{2}
-\tilde\alpha+\frac{n+2}{\tilde q}-\frac{n}{2}=0.
$$
The estimate
\begin{eqnarray}
& &
\biggl\|
|D_x|^{-\tilde\alpha+(n+2)/{\tilde q}-n/2}
\int_0^tU(t-\tau)F(\tau)d\tau
\biggr\|
_{L^\infty({\bf R};L^2({\bf R}^n))}\\
& &
+
\biggl\|
|x|^{-\alpha}
\int_0^tU(t-\tau)F(\tau)d\tau
\biggr\|
_{L^q({\bf R}\times{\bf R}^n)}
\leq
C
\|
|x|^{\tilde\alpha}F
\|_{L^{{\tilde q}'}({\bf R}\times{\bf R}^n)}\nonumber
\end{eqnarray}
holds for radially symmetric $($in $x$$)$ $F$.
}

\vspace{0.3cm}

\noindent{\it Proof of Theorem 2.2}. By Theorem 2.1 we know 
\begin{equation}
\|U^*F\|_{L^2({\bf R}^n)}
\leq
C\||x|^{\tilde\alpha}|D_x|^{-\tilde s}F\|
_{L^{{\tilde q}'}({\bf R}\times{\bf R}^n)}
\end{equation}
and hence
\begin{equation}
\||x|^{-\alpha}|D_x|^s UU^*F\|
_{L^q({\bf R}\times{\bf R}^n)}
\leq
C\||x|^{\tilde\alpha}|D_x|^{-\tilde s}F\|
_{L^{{\tilde q}'}({\bf R}\times{\bf R}^n)}.
\end{equation}
Here
\begin{equation}
(UU^*F)(t,x)
=
U(t)\int_{{\bf R}}U(-\tau)F(\tau)d\tau
=
\int_{{\bf R}}U(t-\tau)F(\tau)d\tau.
\end{equation}
If ${\tilde q}'<q$, 
then the estimate (2.9) is an immediate consequence of 
the Christ-Kiselev lemma \cite{CK} (see also Smith and Sogge \cite{SS}, 
Tao \cite{T}). 
Though the Christ-Kiselev lemma does not apply to the case ${\tilde q}'=q=2$, 
Sugimoto \cite{Su} and Vilela \cite{V} have already generalized the works of 
Kato and Yajima \cite{KY} and Ben-Artzi and Klainerman \cite{BA-K}, 
and they have independently shown the estimate (2.9) 
for ${\tilde q}'=q=2$ without assuming radial symmetry. 
The proof of Theorem 2.2 has been finished.
$\hfill\Box$

\vspace{0.3cm}

\noindent{\it Proof of Theorem 2.3}. 
By the $L^2$ conservation and 
Theorem 2.1 we have 
\begin{eqnarray}
& &
\|UU^*F\|_{L^\infty({\bf R};L^2({\bf R}^n))}\\
& &
=
\|U^*F\|_{L^2({\bf R}^n)}
\leq
C\||x|^{\hat\alpha}|D_x|^{-\hat s}F\|
_{L^{{\hat q}'}({\bf R}\times{\bf R}^n)}.\nonumber
\end{eqnarray}
The estimate (2.10) is an immediate consequence 
of the Christ-Kiselev lemma \cite{CK} and (2.15) as before.
The proof has been finished.
$\hfill\Box$
\section{Proof of Theorem 1.1}
\setcounter{equation}{0}
Our proof starts with the following elementary result.

\vspace{0.3cm}

\noindent{\bf Proposition 3.1.} {\it The inequality 
\begin{equation}
\max
\biggl(
\frac1p,\frac{2}{p-1}-\frac{n-1}{2}
\biggr)
<\frac{2}{p-1}-\frac{n+1}{2p}
\end{equation}
holds for all $n\geq 3$ and 
$1+4/(n+1)<p<1+4/n$.
}

\vspace{0.3cm}

{\it Proof}. 
Though this result is not valid for $n=2$, 
we shall keep $n$ general until the end of our argument 
for clarity. 
By $p_0(n)$ we write the larger root of 
the quadratic equation 
$(n-1)p^2-(n+1)p-2=0$, namely
$$
p_0(n)
=
\frac{n+1+\sqrt{n^2+10n-7}}{2(n-1)}.
$$
By direct computation we can find that 
$1+4/(n+1)<p_0(n)$ for all $n\geq 2$, 
and that $1+4/n<p_0(n)$ for $n=2,3$, 
$1+4/n=2=p_0(n)$ for $n=4$, 
and $p_0(n)<1+4/n$ for $n\geq 5$. 
Because of the equivalence
$$
\frac{1}{p}<\frac{2}{p-1}-\frac{n-1}{2}
\Longleftrightarrow
p<p_0(n),
$$
the proof of Proposition 3.1 is reduced to showing 
\begin{eqnarray}
& &
\frac{2}{p-1}-\frac{n-1}{2}
<
\frac{2}{p-1}-\frac{n+1}{2p},\,\,\,
n=2,3,4,\,\,1+\frac{4}{n+1}<p<1+\frac{4}{n},\\
& &
\frac{2}{p-1}-\frac{n-1}{2}
<
\frac{2}{p-1}-\frac{n+1}{2p},\,\,\,
n\geq 5,\,\,1+\frac{4}{n+1}<p\leq p_0(n),\\
& &
\frac{1}{p}
<
\frac{2}{p-1}-\frac{n+1}{2p},\,\,\,
n\geq 5,\,\,p_0(n)<p<1+\frac4n.
\end{eqnarray}
The inequality (3.2) is equivalent to 
$p>(n+1)/(n-1)$, which 
does not consist with $p<1+4/n$ for $n=2$. 
This is the reason why the case $n=2$ is ruled out in the proposition. 
On the other hand, 
by virtue of $(n+1)/(n-1)\leq 1+4/(n+1)$ for $n\geq 3$, 
the inequality $p>(n+1)/(n-1)$ is automatically satisfied 
if $n\geq 3$. 
For the same reason the statement of (3.3) is true. 
Finally, we see that the inequality (3.4) 
is equivalent to $p<1+4/(n-1)$, 
which is certainly satisfied due to the assumption 
$p<1+4/n$. 
The proof of Proposition 3.1 has been finished.
$\hfill\Box$

\vspace{0.2cm}

In what follows we assume $n\geq 3$. 
Note that 
the assumption $p>1+4/(n+1)$ is equivalent to 
$2/(p-1)-(n-1)/2<1$. 
Therefore, in view of (3.1) we can choose $q_0>2$ so that 
\begin{equation}
\max
\biggl(
\frac{1}{p},
\frac{2}{p-1}-\frac{n-1}{2}
\biggr)
<\frac{2}{q_0}
<
\frac{2}{p-1}-\frac{n+1}{2p},
\end{equation}
and set
\begin{equation}
\alpha_0
:=
\frac{n+2}{q_0}
-
\frac{2}{p-1},\,\,\,
s_0:=
-\alpha_0+\frac{n+2}{q_0}-\frac{n}{2}.
\end{equation}
Note that $s_0$ defined above is equal to 
$-n/2+2/(p-1)$, so that $s_0>0$ for $p<1+4/n$. 
For the purpose of proving Theorem 1.1 
we define the space $X$ of tempered distributions 
on ${\bf R}^{1+n}$ as follows:
\begin{eqnarray}
& &
X:=
\{\,v=v(t,x)\,|\,v\in C({\bf R};{\dot H}^{-s_0}_2({\bf R}^n))
\mbox{\,\,is radially symmetric in\,\,}x,\\
& &
\hspace{3.4cm}
\sup_{t\in{\bf R}}
\||D_x|^{-s_0}v(t,\cdot)\|_{L^2({\bf R}^n)}
+
\||x|^{-\alpha_0}v\|_{L^{q_0}({\bf R}\times{\bf R}^n)}<\infty\,\}.
\nonumber
\end{eqnarray}
We carry out the contraction-mapping argument, using the estimate (2.11) 
with $\tilde q$, $\tilde\alpha$ satisfying 
\begin{equation}
p{\tilde q}'=q_0,\,\,\,-{\tilde\alpha}/p=\alpha_0.
\end{equation}
We rewrite $\tilde q$, $\tilde\alpha$ of (3.8) 
as $q_1$, $\alpha_1$, respectively. 
To proceed, we must verify that 
these satisfy all the conditions of Corollary 2.4. 
It is easy to see that the condition 
$$
-\alpha_0+\frac{n+2}{q_0}-\frac{n}{2}
-\alpha_1+\frac{n+2}{q_1}-\frac{n}{2}=0,
$$
which is reduced to $\alpha_0=(n+2)/q_0-2/(p-1)$, 
is surely satisfied. 
The condition $n/q_0-(n-1)/2<\alpha_0<n/q_0$, 
which is reduced to 
$$
\frac{2}{p-1}-\frac{n-1}{2}<\frac{2}{q_0}<\frac{2}{p-1}
$$
for $\alpha_0=(n+2)/q_0-2/(p-1)$, 
is satisfied by virtue of (3.5). 
The condition $1/q_1\leq 1/2$, 
which is equivalent to $1/p\leq 2/q_0$, 
is assumed in (3.5). 
Since the definition of $q_1$ (see (3.8)) 
implies 
$$
p
\biggl(
\frac{1}{p}-\frac{1}{q_0}
\biggr)
=
\frac{1}{q_1},
$$
the condition $1/q_1>0$ is equivalent to 
$1/q_0<1/p$. 
Because of 
$$
\frac{2}{p-1}-\frac{n+1}{2p}
<
\frac{2}{p}
$$
for all $p>1+4/(n+1)$ 
we surely see that 
the inequality 
$1/q_0<1/p$ is true (see (3.5)). 
The condition $n/q_1-(n-1)/2<\alpha_1$, 
which is equivalent to 
$2/q_0<2/(p-1)-(n+1)/(2p)$, 
is also assumed in (3.5). 
Finally, 
the condition $\alpha_1<n/q_1$, 
which is equivalent to 
$2/(p-1)-n/p<2/q_0$, 
is automatically satisfied 
because $q_0$ has been chosen so that 
$2/(p-1)-(n-1)/2<2/q_0$, 
and we easily see 
$2/(p-1)-n/p<2/(p-1)-(n-1)/2$ 
for $p<1+4/n$. 

Since we have finished the verification of the fact that the exponents 
$q_j$ and $\alpha_j$ $(j=0,1)$ satisfy all the conditions 
in Corollary 2.4, 
let us solve the associated integral equation
\begin{equation}
u(t)
=
U(t)\varphi
-i
\int_0^t
U(t-\tau)f(u(\tau))d\tau
\end{equation}
$(f(u(\tau)):=\lambda|u(\tau)|^{p-1}u(\tau))$ 
by applying the key estimates (2.1) and (2.11). 
An application of (2.1) yields 

\vspace{0.3cm}

\noindent{\bf Proposition 3.2.} {\it 
Suppose $n\geq 3$. 
If $\varphi\in{\dot H}^{-s_0}_2({\bf R}^n)$ 
is radially symmetric, 
then $U(t)\varphi\in X$ and the estimate 
\begin{eqnarray}
& &
\sup_{t\in{\bf R}}
\||D_x|^{-s_0}U(t)\varphi\|_{L^2({\bf R}^n)}
+
\||x|^{-\alpha_0}U(t)\varphi\|
_{L^{q_0}({\bf R}\times{\bf R}^n)}\\
& &
\leq
C\||D_x|^{-s_0}\varphi\|_{L^2({\bf R}^n)}\nonumber
\end{eqnarray}
holds.
}

\vspace{0.3cm}

In view of (3.8) 
we also obtain the following 
by the application of (2.11). 

\vspace{0.3cm}

\noindent{\bf Proposition 3.3.} {\it Suppose $n\geq 3$. 
If $v\in X$, then 
\begin{equation}
\int_0^t U(t-\tau)f(v(\tau))d\tau\in X.
\end{equation}
Moreover, the estimate
\begin{eqnarray}
& &
\sup_{t\in{\bf R}}
\biggl\|
|D_x|^{-s_0}
\int_0^t
U(t-\tau)
\bigl(
f(v_1(\tau))-f(v_2(\tau))
\bigr)d\tau
\biggr\|
_{L^2({\bf R}^n)}\\
& &
\hspace{0.5cm}
+
\biggl\|
|x|^{-\alpha_0}
\int_0^t
U(t-\tau)
\bigl(
f(v_1(\tau))-f(v_2(\tau))
\bigr)d\tau
\biggr\|
_{L^{q_0}({\bf R}\times{\bf R}^n)}\nonumber\\
& &
\leq
C
\bigl(
\||x|^{-\alpha_0}v_1\|_{L^{q_0}({\bf R}\times{\bf R}^n)}
+
\||x|^{-\alpha_0}v_2\|_{L^{q_0}({\bf R}\times{\bf R}^n)}
\bigr)^{p-1}\nonumber\\
& &
\hspace{0.5cm}
\times
\||x|^{-\alpha_0}(v_1-v_2)\|_{L^{q_0}({\bf R}\times{\bf R}^n)}.
\nonumber
\end{eqnarray}
holds for $v_1$,\,$v_2\in X$.
}

\vspace{0.3cm}

\noindent{\it Proof of Proposition 3.3.} Recall that 
\begin{equation}
-\alpha_1+\frac{n+2}{q_1}-\frac{n}{2}
=
\alpha_0-\frac{n+2}{q_0}+\frac{n}{2}
=
-s_0.
\end{equation}
To start with, 
we must note that 
it is actually possible to show 
\begin{equation}
|D_x|^{-s_0}
\int_0^t U(t-\tau)F(\tau)d\tau
\in
\bigl(C\cap L^\infty\bigr)
({\bf R};L^2({\bf R}^n))
\end{equation}
for radially symmetric (in $x$) $F$ 
satisfying 
$|x|^{\alpha_1}F\in L^{q'_1}({\bf R}\times{\bf R}^n)$. 
To verify the continuity in $t$, 
by the density of 
$C_0^\infty({\bf R}\times{\bf R}^n)$ in 
$L^{q'_1}({\bf R}\times{\bf R}^n)$ 
we can choose a sequence 
$\{G_j\}_{j\in{\bf N}}\subset C_0^\infty({\bf R}\times{\bf R}^n)$ 
such that $G_j$ is radially symmetric in $x$ and 
$G_j\to |x|^{\alpha_1}F$ in 
$L^{q'_1}({\bf R}\times{\bf R}^n)$ 
as $j\to\infty$. 
Set $F_j:=|x|^{-\alpha_1}G_j$. 
We claim that 
$|D_x|^{-s_0}F_j\in L^1_{\scriptsize\mbox{loc}}({\bf R};L^2({\bf R}^n))$. 
For the verification of this claim it suffices to show 
$|x|^{-\alpha_1}\in L^q(B_1)$ 
$(n/q=n/2+s_0, 
0<s_0<1/2, 
B_1:=\{x\in{\bf R}^n:|x|<1\})$ 
by the Sobolev embedding theorem 
$L^q({\bf R}^n)\subset {\dot H}^{-s_0}_2({\bf R}^n)$. 
We easily see 
\begin{eqnarray*}
& &
|x|^{-\alpha_1}
\in
L^q(B_1)
\Longleftrightarrow
\alpha_1<\frac{n}{q}
\biggl(
=
\frac{n}{2}+s_0
\biggr)\\
& &
\Longleftrightarrow
\alpha_1-s_0<\frac{n}{2}
\Longleftrightarrow
\frac{n+2}{q_1}-\frac{n}{2}<\frac{n}{2}\,\,\,\,(\mbox{by}\,(3.13))\\
& &
\Longleftrightarrow
1-\frac{p}{q_0}
<
\frac{n}{n+2}\,\,\,\,
\biggl(
\mbox{Recall}\,\,\frac{1}{q_1}=1-\frac{p}{q_0}\,\,\mbox{by}\,(3.8)
\biggr)\\
& &
\Longleftrightarrow
\frac{2}{q_0}>\frac{4}{p(n+2)}.
\end{eqnarray*}
The last inequality obviously holds 
because of the assumptions 
$2/q_0>1/p$ (see (3.5)) and $n\geq 3$. 
Since we have shown 
$|D_x|^{-s_0}F_j\in
L^1_{\scriptsize\mbox{loc}}({\bf R};L^2({\bf R}^n))$, 
we find by the standard argument 
\begin{eqnarray}
& &
|D_x|^{-s_0}
\int_0^t
U(t-\tau)F_j(\tau)d\tau\\
& &
=
\int_0^t
U(t-\tau)|D_x|^{-s_0}F_j(\tau)d\tau
\in
C({\bf R};L^2({\bf R}^n))\nonumber
\end{eqnarray}
$(j=1,2,\dots)$ and hence we conclude by the estimate (2.11) 
that 
this is a Cauchy sequence in 
$\large(C\cap L^\infty\large)({\bf R};L^2({\bf R}^n))$. 
Its limit exists in 
$\large(C\cap L^\infty\large)
({\bf R};L^2({\bf R}^n))$, 
which proves (3.14).

We are in a position to complete the proof of 
Proposition 3.3. 
Since $p q_1'=q_0$ and 
$-\alpha_1/p=\alpha_0$ (see (3.8)), 
the property (3.11) is a direct consequence of 
(2.11), (3.14). 
In view of the basic inequality 
$|f(v_1)-f(v_2)|\leq C(|v_1|+|v_2|)^{p-1}|v_1-v_2|$, 
the estimate (3.12) follows in the same way. 
The proof of Proposition 3.3 has been finished.
$\hfill\Box$

\vspace{0.3cm}

We prove Theorem 1.1. 
For $M>0$ let us define 
$$
X_M:=
\{\,v\in X\,|\,
\sup_{t\in{\bf R}}
\||D_x|^{-s_0}v(t,\cdot)\|_{L^2({\bf R}^n)}
+
\||x|^{-\alpha_0}v\|_{L^{q_0}({\bf R}\times{\bf R}^n)}
\leq M\,\}.
$$
Endowed with the metric
$$
d(v,w):=
\sup_{t\in{\bf R}}
\||D_x|^{-s_0}(v(t,\cdot)-w(t,\cdot))\|_{L^2({\bf R}^n)}
+
\||x|^{-\alpha_0}(v-w)\|_{L^{q_0}({\bf R}\times{\bf R}^n)},
$$
the set $X_M$ is a complete metric space. 
Defining the operator
$$
N(v):=
U(t)\varphi
-i\int_0^t U(t-\tau)f(v(\tau))d\tau\quad(v\in X),
$$
we find by (3.10)--(3.12) that 
there exists a constant $\delta>0$ 
depending on $n$, $p$, and $\lambda$ such that 
if $\||D_x|^{-s_0}\varphi\|_{L^2({\bf R}^n)}\leq\delta$, 
then the operator $N$ has a unique fixed point in $X_{C\delta}$. 
Here $C$ is a suitable positive constant, 
and the fixed point is a solution to the integral equation (3.9) 
which is unique in $X_{C\delta}$. 
We have finished the proof of Theorem 1.1.

\vspace{0.2cm}

{\it Acknowledgements}. The author is grateful to 
Professor Yoshio Tsutsumi for 
showing an interest in this work, 
and to Professor Mitsuru Sugimoto for 
giving a helpful comment on Theorem 1.1. 
He is also grateful for a number of valuable suggestions of the referee. 
He is partly supported by the 
Grant-in-Aid for Young Scientists (B) (No.\,18740069), 
The Ministry of Education, Culture, 
Sports, Science and Technology, Japan.

\begin{flushright}
Kunio Hidano\\
Department of Mathematics\\
Faculty of Education\\
Mie University\\
1577 Kurima-machiya-cho, Tsu\\
Mie Prefecture 514-8507\\
Japan\\
E-mail: hidano@edu.mie-u.ac.jp
\end{flushright}
\end{document}